\documentclass[letterpaper,11pt]{article}
\usepackage[margin=1in]{geometry}  

\usepackage{listings}
\usepackage{amsmath}
\usepackage{amsthm}
\usepackage{tikz}
\usepackage{caption,subcaption}
\usepackage{array}
\usepackage{mdwmath}
\usepackage{multirow}
\usepackage{mdwtab}
\usepackage{eqparbox}
\usepackage{multicol}
\usepackage{amsfonts}
\usepackage{tikz}
\usepackage{multirow,bigstrut,threeparttable}
\usepackage{amsthm}
\usepackage{array}
\usepackage{bbm}
\usepackage{epstopdf}
\usepackage{mdwmath}
\usepackage{mdwtab}
\usepackage{eqparbox}
\usepackage{tikz}
\usepackage{latexsym}
\usepackage{cite}
\usepackage{amssymb}
\usepackage{bm}
\usepackage{amssymb}
\usepackage{graphicx}
\usepackage{mathrsfs}
\usepackage{epsfig}
\usepackage{psfrag}
\usepackage{setspace}
\usepackage[
            CJKbookmarks=true,
            bookmarksnumbered=true,
            bookmarksopen=true,
            colorlinks=true,
            citecolor=blue, 
            linkcolor=blue,
            anchorcolor=blue, 
            urlcolor=blue
            ]{hyperref}
\usepackage{algorithm}
\usepackage{algpseudocode}
\usepackage{stfloats}

\input xy
\xyoption{all}

\theoremstyle{plain}
\newtheorem{theorem}{Theorem}
\newtheorem{lemma}{Lemma}
\newtheorem{assumption}{Assumption}
\newtheorem{corollary}{Corollary}

\theoremstyle{definition}

\def\NN{\mathbb{N}}

\def\calA{\mathcal{A}}

\def\calF{\mathcal{F}}

\def\calL{\mathcal{L}}
\def\calM{\mathcal{M}}

\def\calP{\mathcal{P}}

\def\calS{\mathcal{S}}

\def\calU{\mathcal{U}}

\def\bE{\mathbf{E}}

\def\bP{\mathbf{P}}

\def\bR{\mathbf{R}}

\def\1{\mathbbm{1}}

\usepackage{booktabs}
\usepackage{adjustbox}
\usepackage{multirow}
\usepackage{listings}

\theoremstyle{plain}

\def \bP {\mathbb{P}}
\def \bE {\mathbb{E}}
\def \bR {\mathbb{R}}

\usepackage{xspace}

\newcommand{\stepa}[1]{\overset{\rm (a)}{#1}}
\newcommand{\stepb}[1]{\overset{\rm (b)}{#1}}
\newcommand{\stepc}[1]{\overset{\rm (c)}{#1}}
\newcommand{\stepd}[1]{\overset{\rm (d)}{#1}}
\newcommand{\stepe}[1]{\overset{\rm (e)}{#1}}

\definecolor{myblue}{rgb}{.8, .8, 1}
\definecolor{mathblue}{rgb}{0.2472, 0.24, 0.6} 
\definecolor{mathred}{rgb}{0.6, 0.24, 0.442893}
\definecolor{mathyellow}{rgb}{0.6, 0.547014, 0.24}

\begin{document}

\title{On the High Accuracy Limitation of Adaptive Property Estimation}
\author{Yanjun Han\thanks{Yanjun Han is with the Department of Electrical Engineering, Stanford University, email: \url{yjhan@stanford.edu}. }}

\maketitle

\begin{abstract}
Recent years have witnessed the success of adaptive (or unified) approaches in estimating symmetric properties of discrete distributions, where the learner first obtains a distribution estimator independent of the target property, and then plugs the estimator into the target property as the final estimator. Several such approaches have been proposed and proved to be adaptively optimal, i.e. they achieve the optimal sample complexity for a large class of properties within a low accuracy, especially for a large estimation error $\varepsilon\gg n^{-1/3}$ where $n$ is the sample size. 

In this paper, we characterize the high accuracy limitation, or the penalty for adaptation, for general adaptive approaches. Specifically, we obtain the first known adaptation lower bound that under a mild condition, any adaptive approach cannot achieve the optimal sample complexity for every $1$-Lipschitz property within accuracy $\varepsilon \ll n^{-1/3}$. In particular, this result disproves a conjecture in \cite{acharya2017unified} that the profile maximum likelihood (PML) plug-in approach is optimal in property estimation for all ranges of $\varepsilon$, and confirms a conjecture in \cite{han2021competitive} that their competitive analysis of the PML is tight.
\end{abstract}
\tableofcontents

\section{Introduction and Main Results}
Given $n$ i.i.d. samples drawn from a discrete distribution $p = (p_1,\cdots,p_k)$ of support size $k$, the problem of symmetric (or permutation-invariant) property estimation is to estimate the following quantity
\begin{align*}
	F(p) = \sum_{i=1}^k f(p_i)
\end{align*}
or its variants within a small additive error, for a given function $f: [0,1]\to \bR$. This is a fundamental problem in computer science and statistics with applications in neuroscience \cite{RWDB99}, physics \cite{VBBVP12}, ecology \cite{Chao84, Chao92, BF93, CCGLMCL12}, and others \cite{PW96, ASNRMAS01}.

Over the past decade, there are two main lines of research towards the symmetric property estimation. The first line of research aims to work out the minimax estimation rate and construct the minimax rate-optimal estimators for a given property, including entropy \cite{Paninski2003,Paninski2004,Valiant--Valiant2011,Jiao--Venkat--Han--Weissman2015minimax,wu2016minimax}, support size \cite{Valiant--Valiant2013estimating, wu2019chebyshev}, support coverage \cite{orlitsky2016optimal,ZVVKCSLSDM16}, distance to uniformity \cite{Valiant--Valiant2011power,jiao2018minimax}, sorted $\ell_{1}$ distance \cite{Valiant--Valiant2011power, han2018local}, R\'enyi entropy~\cite{AOST14, AOST17}, nonparametric functionals~\cite{han2020optimal,han2017estimation}, and many others. One of the main findings in these work is that, plugging the empirical distribution into the property often leads to a strictly suboptimal estimator, especially when the function $f$ has some non-smooth parts. They also provided general recipes for the construction of minimax rate-optimal estimators, while the detailed construction crucially depends on the specific property in hand (i.e. classify the smooth and non-smooth parts of $f$, and apply different procedures).

The other line of research aims to achieve a more ambitious goal: find an adaptive (or unified) estimator that achieves the optimal sample complexity for all (or most of) the above symmetric properties. Specifically, the learner aims to obtain a unified distribution estimator $\widehat{p}$ of the true distribution $p$ independent of the property $F$ in hand, and hopes that the plug-in estimator $F(\widehat{p})$ is minimax rate-optimal in estimating $F(p)$ for a large class of properties $F$. This goal may sound too good to be true, for at least two reasons: 
\begin{itemize}
	\item as shown above, the plug-in approach of the empirical distribution, possibly the most natural choice of $\widehat{p}$, does not give the rate-optimal estimator;
	\item the construction of the optimal estimator even for known $F$ is typically quite involved. 
\end{itemize}
However, surprising recent developments show that there does exist such an estimator $\widehat{p}$, and there are even multiple such estimators. One estimator is the \emph{local moment matching} (LMM) estimator in \cite{han2018local} (and its refinement in \cite{han2021competitive}), which is minimax rate-optimal in estimating the true distribution $p$ up to permutation. Moreover, plugging the LMM estimator into the entropy, power sum function, support size, and all $1$-Lipschitz functionals attains the optimal sample complexity for the respective properties within any accuracy $\varepsilon\gg n^{-1/3}$. Another estimator is the \emph{profile maximum likelihood} (PML) estimator proposed in \cite{orlitsky2004modeling}, whose statistical performance was analyzed in \cite{acharya2017unified} via a competitive analysis with an amplification factor $\exp(3\sqrt{n})$ of the error probability; this factor was later improved to $\exp(c'n^{1/3+c})$ for any $c>0$ in \cite{han2021competitive}. Consequently, for a large class of symmetric properties $F$ where there exists a sample-optimal estimator with a sub-Gaussian error probability $\exp(-cn\varepsilon^2)$, the above analyses imply that the PML plug-in approach is also adaptively optimal within any accuracy parameter $\varepsilon \gg n^{-1/3}$. 

These adaptive estimators, albeit promising, still leave some questions. Specifically, we notice the following discrepancy: the estimators constructed in the property-specific manner could achieve the optimal sample complexity for the entire accuracy regime $\varepsilon \gg n^{-1/2}$, while both adaptive estimators above are shown to be optimal only when $\varepsilon\gg n^{-1/3}$. This discrepancy leaves alone the following important question: 

\begin{center}
	\emph{Is there a fundamental limit for general adaptive approaches of property estimation in the high-accuracy regime where $n^{-1/2}\ll \varepsilon\ll n^{-1/3}$?}
\end{center}

Note that there are three possible answers to this question: first, this high-accuracy regime is uncovered simply due to an artifact of the analyses for the above adaptive estimators, and a better theoretical guarantee may be possible. Second, there may exist another fully adaptive estimator which is currently missing. Third, this high-accuracy regime may be a fundamental burden for any adaptive estimator. Specializing this question to the PML, \cite{acharya2017unified} conjectured that ``the PML based approach is indeed optimal for all ranges of $\varepsilon$", while \cite{han2021competitive} conjectured that $\varepsilon\gg n^{-1/3}$ is the best possible range for the PML to be adaptively optimal. However, even for the PML, which is a specific choice of an adaptive estimator, the lower bound analysis is missing. 

In this paper, we show that the latter conjecture is true even for general adaptive estimation: there is a phase transition for the performance of adaptive estimators at the accuracy parameter $\varepsilon\asymp n^{-1/3}$, while beyond this point, there is an unavoidable price that \emph{any} adaptive estimator needs to pay on the sample complexity. In other words, for a reasonable family of symmetric properties, although property-specific approaches are optimal for the full accuracy range $\varepsilon\gg n^{-1/2}$, any adaptive approach fails to achieve the optimal sample complexity for at least \emph{one} of the properties if $\varepsilon\ll n^{-1/3}$. Specifically, our main contributions are as follows: 

\begin{enumerate}
	\item We prove the first tight adaptation lower bound for the class of all $1$-Lipschitz properties. We show that although the sample complexity for each $1$-Lipschitz property is at most $O(k/(\varepsilon^2\log k))$ for any $\varepsilon\gg n^{-1/2}$, under a mild assumption, any adaptive estimator must incur a sample complexity at least $\Omega(k/\varepsilon^2)$ for every $\varepsilon\ll n^{-1/3}$. 
	\item As a corollary, we obtain a tight competitive analysis for the PML plug-in approach. Specifically, we show that the amplification factor of the error probability in the PML competitive analysis is at least $\exp(\Omega(n^{1/3-c}))$ for every $c>0$, resolving the tightness conjecture of the upper bound $\exp(O(n^{1/3+c}))$ for every $c>0$ in \cite{han2021competitive}. 
	\item We consider a new class of adaptive estimation problems, where we aim to adapt to a family of loss functions instead of the parameter sets in the traditional setting. We propose a generalized Fano's inequality to establish the adaptation lower bound for the new problem, which could be of independent interest. 
\end{enumerate}

The rest of the paper is organized as follows. We introduce the necessary notations in Section \ref{subsec.notation}, and state the main adaptation lower bound in Section \ref{subsec.result}. We also compare our setting and results with an extensive set of prior work in Section \ref{subsec.relatedwork}. For the proofs of main results, Section \ref{sec.theorem1} and \ref{sec.theorem2} present the detailed proofs of Theorems \ref{thm:main} and \ref{thm:PML}, respectively. In particular, Section \ref{subsec.highlevel} presents a novel and general idea in proving the adaptation lower bound, and Section \ref{subsec.fano} presents a generalized Fano's inequality. Conclusions and open problems are drawn in Section \ref{sec.conclusion}. 

\subsection{Notations}\label{subsec.notation}
Throughout the paper we adopt the following notations. Let $\NN$ be the set of all positive integers, and for $n\in \NN$, let $[n]\triangleq \{1,2,\ldots,n\}$. For a finite set $A$, let $|A|$ be the cardinality of $A$. For $k\in\NN$, let $\calM_k$ be the set of all discrete distributions supported on $[k]$. For two probability measures $P,Q$ on the same probability space, let $\|P-Q\|_{\text{TV}} = \int |dP-dQ|/2$, $D_{\text{KL}}(P \| Q) = \int dP\log(dP/dQ)$, and $\chi^2(P\|Q) = \int (dP - dQ)^2/dQ$ be the total variation (TV) distance, the Kullback--Leibler (KL) divergence, and the $\chi^2$-divergence between $P$ and $Q$, respectively. For random variables $X$ and $Y$, let $I(X;Y) = \int dP_{XY}\log(dP_{XY}/dP_X\otimes dP_Y)$ be the mutual information. For $p\in \calM_k$, let $\bP_p$ and $\bE_p$ denote the probability and expectation taken with respect to the i.i.d. samples $X_1,\cdots,X_n\sim p$, respectively. We also adopt the following asymptotic notations. For non-negative sequences $\{a_n\}$ and $\{b_n\}$, we write $a_n\lesssim b_n$ (or $a_n = O(b_n)$) to denote that $\limsup_{n\to\infty}a_n/b_n<\infty$, and $a_n \gtrsim b_n$ (or $a_n = \Omega(b_n)$) to denote $b_n\lesssim a_n$, and $a_n\asymp b_n$ (or $a_n = \Theta(b_n)$) to denote both $a_n\lesssim b_n$ and $b_n\lesssim a_n$. We also write $a_n\ll b_n$ to denote that $\limsup_{\varepsilon\to 0}\limsup_{n\to\infty} n^{\varepsilon} a_n/b_n = 0$, and $a_n \gg b_n$ to denote $b_n\ll a_n$. 

\subsection{Main Results}\label{subsec.result}
To state our main adaptation lower bound, we first need to define the family of adaptive estimators as well as the family of symmetric property estimation problems. In this paper, we consider the class $\calF_{\text{Lip}}$ of all $1$-Lipschitz properties expressed as $F(p) = \sum_{i=1}^k f(p_i)$ with a $1$-Lipschitz function $f: [0,1] \to \mathbb{R}$, i.e. $|f(x) - f(y)|\le |x-y|$ for all $x,y\in [0,1]$. As for the class of adaptive estimators, we require that the learner obtains a single discrete distribution estimator $\widehat{p} = (\widehat{p}_1,\cdots,\widehat{p}_k)$ based on the observations $X^n$, and then uses the plug-in estimator $F(\widehat{p})$ to estimate the property $F(p)$. To measure the performance of this adaptive estimator, we consider the expected estimation error $\bE_p |F(\widehat{p}) - F(p)|$ for the worst-case discrete distribution $p\in \calM_k$ \emph{and} the worst-case $1$-Lipschitz property $F\in \calF_{\text{Lip}}$. In other words, in this paper we are interested in characterizing the following \emph{adaptive minimax risk}: 
\begin{align}\label{eq:adaptive_minimax_risk}
R_{\text{adaptive}}^\star(n,k) \triangleq \inf_{\widehat{p}}\sup_{F \in \calF_{\text{Lip}}}\sup_{p\in \calM_k} \bE_{p} |F(\widehat{p}) - F(p)|.
\end{align}

For technical reasons, we also assume the following mild assumption for the single distribution estimator $\widehat{p}$. 
\begin{assumption}\label{assump}
For each $n,k\in \NN$, we assume that the distribution estimator $\widehat{p}(X^n) = (\widehat{p}_1,\cdots,\widehat{p}_k)$ satisfies (where $\calS_k$ denotes the permutation group over $[k]$)
\begin{align*}
\sup_{p\in \calM_k} \bE_p \left[\min_{\sigma\in \calS_k} \sum_{i=1}^k |\widehat{p}_{\sigma(i)} - p_i| \right] \le A(n)\cdot \sqrt{\frac{k}{n}}, 
\end{align*}
with $A(n) \ll n^{\delta}$ for every $\delta > 0$. We will use $\calP$ to denote the class of all such estimators $\widehat{p}$. 
\end{assumption}

Assumption \ref{assump} essentially requires that the single distribution estimator $\widehat{p}$ used in the adaptive approach must be a \emph{reasonably} good estimator of the true distribution $p$ up to permutation, where the term \emph{reasonably} means that the estimator cannot be much worse than the empirical estimator. We provide three reasons on why we believe this assumption to be mild. First, it is very natural to expect or require that a good distribution estimator used in the adaptive approach should be sound not only after being plugging into various properties, but also \emph{before} the plug-in process in terms of the (sorted) distribution estimation. In other words, Assumption \ref{assump} could be treated as an additional requirement for any sound adaptive approach. Second, Assumption \ref{assump} holds for many natural or known estimators. For example, the empirical distribution satisfies Assumption \ref{assump} with $A(n)\equiv 1$ (see, e.g. \cite{han2015minimax}), and both known adaptive estimators, i.e. LMM and PML, also belong to $\calP$ with $A(n) = \text{polylog}(n)$ (cf. \cite{han2018local} for the LMM, and the proof of Theorem \ref{thm:PML} for the PML). Hence, restricting to the estimator class $\calP$ still leads to non-trivial lower bounds for these known estimators. Third, a larger quantity $A(n)$ in Assumption \ref{assump} only shrinks the accuracy regime from $\varepsilon \ll n^{-1/3}$ to $\varepsilon\ll (nA(n))^{-1/3}$, but does not affect the claimed minimax lower bound in the new accuracy regime. In addition to these reasons, we remark that Assumption \ref{assump} is mostly a technical assumption, and we conjecture that the following Theorem \ref{thm:main} still holds without it. 

Restricting to the estimator class $\calP$, the following theorem characterizes the tight adaptive minimax rate for $1$-Lipschitz property estimation. 
\begin{theorem}\label{thm:main}
For each $n,k\in \NN$, it holds that 
\begin{align*}
\inf_{\widehat{p}\in \calP}\sup_{F \in \calF_{\text{\rm Lip}}}\sup_{p\in \calM_k} \bE_{p} |F(\widehat{p}) - F(p)| \asymp \begin{cases}
\sqrt{\frac{k}{n\log n}} & \text{if } n^{1/3}\ll k\lesssim n\log n, \\
\sqrt{\frac{k}{n}} & \text{if } 1 \ll k\ll n^{1/3}.
\end{cases}
\end{align*}
\end{theorem}

Theorem \ref{thm:main} can also be equivalently formulated in terms of the optimal sample complexity. 
\begin{corollary}\label{cor.sample_complexity}
It is sufficient and necessary to have $n=\Theta(k/(\varepsilon^2\log k))$ samples for the existence of an adaptive estimator in $\calP$ to estimate all $1$-Lipschitz properties within error $\varepsilon$ if $\varepsilon \gg n^{-1/3}$, and it is sufficient and necessary to have $n=\Theta(k/\varepsilon^2)$ samples for the existence of an adaptive estimator in $\calP$ to estimate all $1$-Lipschitz properties within error $\varepsilon$ if $n^{-1/2}\ll \varepsilon \ll n^{-1/3}$. 
\end{corollary}

Let us appreciate the result of Theorem \ref{thm:main} via the comparison with other results. First, there will be no phase transition in the high-accuracy regime if we do \emph{not} require an adaptive estimator. Specifically, the following result was shown in \cite{hao2019unified}:\footnote{The original paper did not require to use a plug-in estimator, but any estimator $\widehat{F}$ could be clipped to the range of $F(\cdot)$ and written as $F(\widehat{p})$ for some $\widehat{p}\in \calM_k$.} 
\begin{align}\label{eq:non-adaptive}
\sup_{F \in \calF_{\text{\rm Lip}}}\inf_{\widehat{F}}\sup_{p\in \calM_k} \bE_{p} |\widehat{F} - F(p)| \asymp \sqrt{\frac{k}{n\log n}}, \quad 1\ll k\lesssim n\log n.
\end{align}
Comparing Theorem \ref{thm:main} and \eqref{eq:non-adaptive}, we observe that simply after a swap of the infimum and supremum, the minimax risk becomes significantly different. In particular, there is a \emph{strict separation} between the best achievable errors for adaptive and non-adaptive approaches, and the learner may need to pay a strict penalty on the estimation error to achieve adaptation. 

Second, we also compare Theorem \ref{thm:main} with a similar form of the minimax risk in the problem of estimating sorted distribution, where \cite{han2018local} shows that
\begin{align}\label{eq:est_sorted_dist}
\inf_{\widehat{p}}\sup_{p\in \calM_k} \bE_{p} \left[ \sup_{F \in \calF_{\text{\rm Lip}}}|F(\widehat{p}) - F(p)| \right] \asymp \begin{cases}
\sqrt{\frac{k}{n\log n}} & \text{if } n^{1/3}\ll k\lesssim n\log n, \\
\sqrt{\frac{k}{n}} & \text{if } 1 \ll k\ll n^{1/3}.
\end{cases}
\end{align}
As $\bE[\sup_n X_n]\ge \sup_n\bE[X_n]$, the quantity in \eqref{eq:est_sorted_dist} is no smaller than our adaptive minimax risk in \eqref{eq:adaptive_minimax_risk}, and thus implies the upper bound in Theorem \ref{thm:main}. However, the lower bound of Theorem \ref{thm:main} is the most challenging part and stronger than what \eqref{eq:est_sorted_dist} gives. In particular, we remark that after exchanging the expectation and supremum, the lower bound argument will become fundamentally different, and the traditional approaches fail to give the tight adaptive lower bound. We refer to Section \ref{subsec.highlevel} for more details. Moreover, comparing the results of Theorem \ref{thm:main} and \eqref{eq:est_sorted_dist}, we remark that \eqref{eq:est_sorted_dist} gives a tight phase transition for the \emph{problem}, while the adaptive lower bound in Theorem \ref{thm:main} shows a tight phase transition \emph{only for adaptive approaches}. Technically, the former transition could be derived by studying different regimes of the problem, while the latter transition requires to also take into account the crucial nature of the adaptive approach. 

The general adaptive lower bound of Theorem \ref{thm:main} also gives tight and non-trivial lower bounds for some known adaptive approaches. For example, for the LMM adaptive approach in \cite{han2018local}, Theorem \ref{thm:main} shows that the condition $\varepsilon\gg n^{-1/3}$ required for its optimality in property estimation is not superfluous, but in general unavoidable. The implication for the PML adaptive approach \cite{orlitsky2004modeling} is even more surprising; to fully describe this we need to recall some basics of PML.

Given $n$ i.i.d. observations $X_1,\cdots,X_n$ drawn from a discrete distribution $p$ supported on the domain $[k]$, the \emph{profile} of the observations is defined as a vector $\phi = (\phi_0,\cdots,\phi_n)$ with $\phi_i$ being the number of domain elements $j\in [k]$ which appear exactly $i$ times in the sample. For example, $\phi_0$ is the number of unseen elements, and $\phi_1$ is the number of unique elements, i.e. appearing exactly once. Let $\Phi_{n,k}$ be the set of all possible profiles with $n$ observations and support size $k$. Note that for any $\phi\in \Phi_{n,k}$ and $p\in \calM_k$, we could compute the probability that the resulting profile is $\phi$ under true distribution $p$, denoted by $\bP(p,\phi)$. The \emph{profile maximum likelihood} (PML) distribution estimator is then defined as
\begin{align*}
p^{\text{PML}}(\phi) = \arg\max_{p\in \calM_k} \bP(p, \phi).
\end{align*}
In other words, upon observing the profile $\phi$, the PML estimator is the discrete distribution which maximizes the probability of observing $\phi$. This estimator is interesting in several aspects. From the optimization side, the probability $\bP(p,\phi)$ is a highly non-convex function of $p$, and it is very challenging to compute the exact or approximate PMLs. From the statistical side, as $\bP(p,\phi)$ does not admit an additive form even under i.i.d. models (unlike the traditional log-likelihood), even first-order asymptotic properties are challenging to establish for the PML. After 13 years of its invention, a useful statistical property of the PML was established in \cite{acharya2017unified} in terms of an interesting \emph{competitive analysis}: for every property $F$ and accuracy parameter $\varepsilon$, it holds that
\begin{align}\label{eq:competitive}
\sup_{p\in \calM_k} \bP_p( |F(p^{\text{PML}})- F(p)| \ge 2\varepsilon )\le \exp(3\sqrt{n})\cdot \inf_{\widehat{F}}\sup_{p\in \calM_k} \bP_p( |\widehat{F} - F(p)| \ge \varepsilon ). 
\end{align}
Specifically, \eqref{eq:competitive} gives an indirect statistical analysis of the PML plug-in approach which depends on the performance of another estimator. For many properties (such as all $1$-Lipschitz properties), the minimax error probability in the RHS of \eqref{eq:competitive} behaves as $\exp(-\Omega(n\varepsilon^2))$ when $n$ exceeds the optimal sample complexity, thus \eqref{eq:competitive} shows that the PML plug-in approach is adaptively optimal for $\varepsilon\gg n^{-1/4}$. The proof of \eqref{eq:competitive} used only the defining property of PML in a delicate way, and the error amplification factor $\exp(3\sqrt{n})$ follows from a simple union bound and a cardinality bound on the number of profiles $|\Phi_{n,k}|\le \exp(3\sqrt{n})$. 

The paper \cite{acharya2017unified} asked whether the above error amplification factor $\exp(3\sqrt{n})$ could be improved in general; three years later \cite{han2021competitive} provided an affirmative answer. Specifically, using a chaining property of the PML distributions, \cite{han2021competitive} showed the following improvement
\begin{align}\label{eq:competitive_improved}
\sup_{p\in \calM_k} \bP_p( |F(p^{\text{PML}})- F(p)| \ge (2+o(1))\varepsilon )\le \exp(c'n^{1/3+c})\cdot \left(\inf_{\widehat{F}}\sup_{p\in \calM_k} \bP_p( |\widehat{F} - F(p)| \ge \varepsilon )\right)^{1-c}
\end{align}
for any absolute constant $c>0$ and some $c'>0$ depending only on $c$. Using \eqref{eq:competitive_improved}, the accuracy range for the optimality of PML could be improved to $\varepsilon\gg n^{-1/3}$ for the aforementioned properties. It was also conjectured in \cite{han2021competitive} that the new amplification factor in \eqref{eq:competitive_improved} is tight, but little intuition was provided. 

Surprisingly, without directly analyzing the PML adaptive approach, Theorem \ref{thm:main} implies the tightness of the error amplification factor in \eqref{eq:competitive_improved}, as summarized in our next main theorem. 

\begin{theorem}\label{thm:PML}
For any given constants $C>0, c_1\in (0,1/3)$ and $c_2\in (0,1)$, it holds that 
\begin{align*}
\liminf_{n\to\infty} n^{-(1/3-c_1)}\cdot \sup_{F\in \calF_{\text{\rm Lip}}}\sup_{k,\varepsilon>0}\log \frac{\sup_{p\in \calM_k} \bP_p( |F(p^{\text{\rm PML}})- F(p)| \ge C\varepsilon )}{\left(\inf_{\widehat{F}}\sup_{p\in \calM_k} \bP_p( |\widehat{F} - F(p)| \ge \varepsilon )\right)^{1-c_2}} = +\infty. 
\end{align*}
\end{theorem}

After some algebra, it is clear that Theorem \ref{thm:PML} rules out the possibility that the exponent $O(n^{1/3+c})$ of the amplification factor in \eqref{eq:competitive_improved} could be improved to $O(n^{1/3-c})$ in general. Therefore, Theorem \ref{thm:PML} implies that the general competitive analysis of the PML in \cite{han2021competitive} is essentially tight, thereby resolves the conjecture in \cite{han2021competitive}. 

We provide two additional remarks on Theorem \ref{thm:PML}. First, the validity of Theorem \ref{thm:PML} is irrelevant to Assumption \ref{assump}, as the PML estimator is a simple instance which satisfies Assumption \ref{assump}. Second, the lower bound in Theorem \ref{thm:PML} does \emph{not} rule out the possibility that the PML adaptive approach could be fully optimal for \emph{some} property. For example, it was shown in \cite{charikar2019general} that the PML plug-in approach is fully optimal in estimating the support size. It will be an understanding open question to propose a tight analysis of the PML estimator for specific properties. 

\subsection{Related Work}\label{subsec.relatedwork}
\textbf{Property Estimation.} There has been a rich line of research towards the optimal estimation of properties (or functionals) of high-dimensional parameters, especially in the past decade. Starting from some early work \cite{Lepski--Nemirovski--Spokoiny1999estimation,Paninski2003,Paninski2004,Cai--Low2011,Valiant--Valiant2011,Valiant--Valiant2011power,Valiant--Valiant2013estimating}, the fully minimax rate-optimal estimators in all accuracy regimes were obtained for the Shannon entropy in \cite{Jiao--Venkat--Han--Weissman2015minimax,wu2016minimax}. They also provided general recipes for both the estimator construction and tight minimax lower bounds. Specifically, the crux of the optimal estimator construction lies in the classification of smooth and non-smooth regimes and the usage of polynomial approximation to reduce bias in the non-smooth regime, and the minimax lower bound relies on the duality between moment matching and best polynomial approximation. Since then, these general recipes together with their non-trivial extensions have been applied to various other properties, e.g. the R\'{e}nyi entropy \cite{AOST14,AOST17}, support size \cite{wu2019chebyshev}, support coverage \cite{orlitsky2016optimal,ZVVKCSLSDM16,polyanskiy2019dualizing}, distance to uniformity \cite{jiao2018minimax}, general $1$-Lipschitz property \cite{hao2019broad,hao2019unified}, $L_1$ distance \cite{jiao2018minimax}, KL divergence \cite{bu2018estimation,han2020minimax}, and nonparametric functionals \cite{han2017estimation,han2020optimal}. We refer to the survey \cite{verdu2019empirical} for an overview of these results. There is also another line of recent work on estimating a population of parameters or distribution under a Wasserstein distance, a problem closely related to property estimation, via projection-based methods without explicit polynomial approximation \cite{kong2017spectrum,tian2017learning,han2018local,rigollet2019uncoupled,vinayak2019maximum,vinayak2019optimal,wu2020optimal,jana2020extrapolating}. While the above work completely characterized the complexity of many \emph{given problems} in property estimation, the complexity of adaptive estimation in \emph{a set of such problems} is largely missing. For example, the $\Omega(\sqrt{k/(n\log n)})$ lower bound for large $k$ in Theorem \ref{thm:main} simply follows from the complexity of estimating a particular $1$-Lipschitz property, but the main $\Omega(\sqrt{k/n})$ lower bound for small $k$ becomes the crucial complexity of adaptive approaches and thus does not follow from the above set of results or tools. 
\vspace{1em}

\noindent\textbf{Adaptive Property Estimation.} More recently the problem of adaptive, or unified, property estimation has drawn several research attention. As reviewed in the introduction, possibly the most well-known adaptive approach is the PML plug-in approach, with early statistical developments in \cite{orlitsky2004modeling,orlitsky2011estimating,anevski2017estimating}. Since \cite{acharya2017unified} provided the first competitive analysis of the PML plug-in approach, there have been several follow-up papers on the statistical analysis of the PML. Some work focused on the application of the competitive analysis and the construction of the estimator achieving the minimax error probability in \eqref{eq:competitive}, e.g. \cite{hao2019broad}. Some work focused on proper modifications of the PML to achieve better adaptation, e.g. \cite{hao2019broad,charikar2019general}; however, these modified distribution estimators will depend on the target property and are thus not fully unified. Other work aimed to improve the competitive analysis in \cite{acharya2017unified}; for example, \cite{hao2020profile} obtained a distribution-dependent amplification factor without changing the worst-case analysis, and \cite{han2021competitive} improved this factor to $\exp(O(n^{1/3+c}))$ in general. However, none of the above work studied the limitation of the PML plug-in approach, even for concrete examples. Therefore, the lower bound analysis, especially the possible separation compared with the optimal estimator, of the PML is missing. 

Another adaptive approach plugs in the LMM estimator proposed in \cite{han2018local}. Different from the general competitive analysis of PML, the performance of the LMM approach could be directly analyzed for given properties based on its moment matching performance in each local interval. Built on the LMM performance analysis in estimating entropy, power sum function, and support size, the authors of \cite{han2018local} commented that the LMM pays some penalty for being a unified approach. However, this comment was only an insight, and there was no lower bound to support it rigorously. The current work fills in this gap and shows that the price observed for the LMM is in fact unavoidable even for general adaptive approaches. 
\vspace{1em}

\noindent\textbf{Adaptation Lower Bound.} We also review and compare with some known tools to establish adaptation lower bounds, mainly taken from the statistics literature. Adaptation is an important topic in statistics; for example, in nonparametric estimation one may aim to design a density estimator adapting to different smoothness parameters, or in hypothesis testing one may wish to propose an adaptive test procedure against several different alternatives. However, for some problems the adaptation could be achieved without paying any penalty (e.g. density estimation \cite{lepskii1992asymptotically,donoho1995wavelet}, $L_r$ norm estimation with non-even $r$ \cite{han2017estimation}), while for others some adaptation penalties are inevitable (e.g. linear \cite{efromovich1994adaptive} or quadratic \cite{efromovich1996optimal} functional of densities). The main technical tool to establish tight penalties of adaptation is the \emph{constrained risk inequality} originally developed in \cite{brown1996constrained} and generalized in \cite{Cai--Low2011,duchi2018constrained}. Roughly speaking, this type of inequality asserts that if an estimator achieves a too small error at one point, it must incur a too large error at another point; therefore, adaptation may incur a penalty as it might be required to adapt to easier problems and achieve a too small error. For testing, there is also another approach to establish adaptation lower bounds, where the key is to use a mixture of different alternative distributions which could be closer to the null than any fixed alternative; see \cite{spokoiny1996adaptive} and also \cite[Chapter 8]{gine2016mathematical} for examples. 

However, we remark that our adaptive estimation problem is fundamentally different. In the above work, the target of adaptive estimation is to adapt to different (usually a nest of) \emph{parameter sets}, such as H\"{o}lder balls with different smoothness parameters. In contrast, we consider a fixed parameter set (i.e. $p\in \calM_k$), but wish to adapt to different \emph{loss functions} for the final estimator. Establishing adaptation lower bounds for different losses is novel to our knowledge, and the above tools are not applicable in this problem. Consequently, we aim to provide useful tools (cf. Section \ref{subsec.highlevel}) for this new adaptation problem, and expect them to be a helpful addition to the literature on adaptive estimation. 

\section{Proof of Theorem \ref{thm:main}}\label{sec.theorem1}
This section is devoted to the proof of Theorem \ref{thm:main}. Note that the upper bound is achieved by the LMM estimator for $k\gg n^{1/3}$ and the empirical distribution for $k\ll n^{1/3}$ \cite{han2018local}\footnote{Note that \cite{han2018local} shows that both the LMM and empirical distributions belong to $\calP$.}, and the lower bound for $k\gg n^{1/3}$ follows from the minimax lower bound for estimating a specific $1$-Lipschitz property, i.e. the distance to uniformity $F(p) = \sum_{i=1}^k |p_i - 1/k|$ \cite{jiao2018minimax}. Therefore, it remains to prove the following adaptation lower bound:
\begin{align}\label{eq:thm1_target}
\inf_{\widehat{p}\in \calP}\sup_{F \in \calF_{\text{\rm Lip}}}\sup_{p\in \calM_k} \bE_{p} |F(\widehat{p}) - F(p)| \gtrsim \sqrt{\frac{k}{n}}, \quad 1\ll k\ll n^{1/3}. 
\end{align}

This section is organized as follows. Section \ref{subsec.highlevel} presents a high-level overview of the idea in proving the adaptation lower bounds to a class of loss functions in an abstract decision-theoretic setup, and Section \ref{subsec.fano} introduces a generalized Fano's inequality for the adaptation lower bounds. The details to feed into these tools are worked out in Section \ref{subsec.detailedproof}. 

\subsection{High-level Idea}\label{subsec.highlevel}
We consider a general decision-theoretic setup \cite{Wald1950statistical}. Let $(P_\theta)_{\theta\in\Theta}$ be a general statistical model with parameter set $\Theta$, and $\calA$ be the space of all possible actions the learner could take. In other words, the learner obtains an observation $X\sim P_\theta$ with some unknown $\theta$, and then maps $X$ to a random action $a(X)\in \calA$. Let $L: \Theta\times \calA \to \bR_+$ be any (measurable) loss function, the problem of \emph{minimax estimation} is to characterize the following minimax risk:
\begin{align}\label{eq:minimax}
R^\star(\Theta,\calA,L) = \inf_{a}\sup_{\theta\in\Theta} \bE_\theta[L(\theta, a(X))]. 
\end{align}
Similarly, the problem of \emph{adaptive minimax estimation} with respect to a class of loss functions $\calL$ is to characterize the following adaptive minimax risk: 
\begin{align}\label{eq:minimax_adaptive}
R^\star(\Theta,\calA,\calL) = \inf_{a}\sup_{\theta\in\Theta}\sup_{L\in\calL} \bE_\theta[L(\theta, a(X))]. 
\end{align}
To see how \eqref{eq:minimax_adaptive} is related to the adaptive property estimation, we could set $P_\theta$ to be the distribution of $n$ i.i.d. samples from the discrete distribution $\theta$, with $\Theta = \calM_k$. Moreover, $\calA = \calM_k$, $L_F(\theta,a) = |F(\theta) - F(a)|$, and $\calL = \{L_F: F \text{ is a 1-Lipschitz property} \}$.

There are several well-known tools to establish the lower bound of \eqref{eq:minimax}, where a standard and prominent tool is the reduction to hypothesis testing problems; see, e.g. \cite{yu1997assouad,Tsybakov2009}. The main step is to find $\theta_1,\cdots,\theta_M\in \Theta$ such that both the \emph{separation condition} and the \emph{indistinguishability condition} hold: the separation condition typically requires that $\inf_{a\in \calA} [L(\theta_i, a) + L(\theta_j, a)]\ge \Delta$ for some separation parameter $\Delta>0$ and all $i\neq j$, and the indistinguishability condition essentially states that any learner could not determine the true parameter $\theta_i$ based on her observations if the truth $i\in [M]$ is chosen uniformly at random. Then it might be tempting to think that one only needs to replace $L(\theta,a)$ by $\sup_{L\in \calL} L(\theta,a)$ in the above arguments to lower bound \eqref{eq:minimax_adaptive}. However, this approach will place the supremum in $L$ inside the expectation in \eqref{eq:minimax_adaptive}, and thus provide a lower bound for a larger quantity like \eqref{eq:est_sorted_dist}. An alternative way is to use the trivial inequality $R^\star(\Theta,\calA,\calL) \ge \sup_{L\in \calL}R^\star(\Theta,\calA,L)$ and then lower bound the latter quantity. Although this gives a valid lower bound, it is not strong enough in our problem where $R^\star(\Theta,\calA,\calL) \gg\sup_{L\in \calL}R^\star(\Theta,\calA,L)$ in view of Theorem \ref{thm:main} and \eqref{eq:non-adaptive}. 

The main idea to fix the above difficulty is that in addition to choose $M$ points $\theta_1,\cdots,\theta_M \in \Theta$ corresponding to different statistical models, we also find $M$ different loss functions $L_1,\cdots,L_M\in \calL$ tailored for the respective models. Specifically, the indistinguishability condition is unchanged as it depends only on $\theta_1,\cdots,\theta_M$, while the separation condition could be replaced by $\inf_{a\in \calA} [L_i(\theta_i, a) + L_j(\theta_j, a)]\ge \Delta$ for all $i\neq j$. Despite its simplicity, this idea gives the tight adaptation lower bound for the property estimation, and can thus be viewed as an adaptive version of the hypothesis testing approach for the adaptation lower bound. 

\subsection{Generalized Fano's Inequality}\label{subsec.fano}
There is an additional difficulty to apply the aforementioned high-level idea to our problem, i.e. the new separation condition $L_i(\theta_i, a) + L_j(\theta_j, a) \ge \Delta$ does not hold for any action $a\in \calA$, but instead holds for the random action $a(X)$ with a strictly positive probability. To account for this subtlety, we propose the following version of the Fano's inequality. 

\begin{lemma}[Generalized Fano's Inequality]\label{lemma.fano_generalized}
In the above decision-theoretic setup, suppose that $\theta_1,\cdots,\theta_M\in \Theta$ and $L_1,\cdots,L_M\in \calL$ are chosen. Assume that there exists $\calA_0 \subseteq \calA$ such that
\begin{align*}
\inf_{a\in \calA_0} \left[ L_i(\theta_i, a) + L_j(\theta_j, a) \right] \ge \Delta > 0, \quad \forall i,j \in [M], i\neq j, 
\end{align*}
and an estimator $a(X)$ satisfies that $P_{\theta_i}(a(X) \in \calA_0) \ge p_{\min} > 0$ for all $i\in [M]$. Then for this estimator we have
\begin{align*}
\sup_{\theta\in\Theta}\sup_{L\in\calL} \bE_\theta[L(\theta,a(X))] \ge \frac{\Delta}{2}\left(p_{\min} - \frac{I(U;X)+p_{\min}\log 2}{\log M}\right), 
\end{align*}
where $I(U;X)$ denotes the mutual information between $U\sim \text{\rm Uniform}([M])$ and $X\mid U\sim P_{\theta_U}$. 
\end{lemma}

Note that when $L_i\equiv L$ and $p_{\min} = 1$, Lemma \ref{lemma.fano_generalized} reduces to the traditional Fano's inequality \cite{Cover--Thomas2006}. Hence, Lemma \ref{lemma.fano_generalized} is a generalization of the Fano's inequality in the sense that it gives an adaptation lower bound with a soft separation condition. We prove Lemma \ref{lemma.fano_generalized} in the remainder of this subsection. First, as the maximum is no smaller than the average, we have
\begin{align}\label{eq:max_to_avg}
\sup_{\theta\in\Theta}\sup_{L\in\calL} \bE_\theta[L(\theta,a(X))] \ge \frac{1}{M}\sum_{i=1}^M \bE_{\theta_i}[L_i(\theta_i,a(X))]. 
\end{align}
For each $i\in [M]$, let $Q_{i}$ be the conditional distribution of $a(X)$ with $X\sim P_{\theta_i}$ conditioning on the event $a(X) \in \calA_0$. Then by the non-negativity of each $L_i$ and definition of $p_{\min}$, 
\begin{align*}
\bE_{\theta_i}[L_i(\theta_i,a(X))] \ge P_{\theta_i}(a(X)\in \calA_0) \cdot \bE_{a\sim Q_i}[L_i(\theta_i, a)]\ge p_{\min}\cdot \bE_{a\sim Q_i}[L_i(\theta_i, a)],
\end{align*}
and therefore \eqref{eq:max_to_avg} gives
\begin{align}\label{eq:conditional}
\sup_{\theta\in\Theta}\sup_{L\in\calL} \bE_\theta[L(\theta,a(X))] \ge p_{\min}\cdot \frac{1}{M}\sum_{i=1}^M \bE_{a\sim Q_i}[L_i(\theta_i, a)]. 
\end{align}
The next few steps are similar to the proof of the traditional Fano's inequality. For each $a\in \calA_0$, define a test $\Psi(a) = \arg\min_{i \in [M]}L_i(\theta_i,a)$. Then by the separation condition, we have
\begin{align*}
L_i(\theta_i, a) \ge \frac{L_i(\theta_i, a) + L_{\Psi(a)}(\theta_{\Psi(a)}, a) }{2} \ge \frac{\Delta}{2}\cdot \1(\Psi(a)\neq i), \quad \forall i\in [M], a\in \calA_0,
\end{align*}
and therefore \eqref{eq:conditional} gives
\begin{align}\label{eq:test}
\sup_{\theta\in\Theta}\sup_{L\in\calL} \bE_\theta[L(\theta,a(X))] \ge \frac{\Delta p_{\min}}{2}\cdot \frac{1}{M}\sum_{i=1}^M Q_i(\Psi(a) \neq i) \ge \frac{\Delta p_{\min}}{2}\left(1 - \frac{I(U;Y)+\log 2}{\log M}\right), 
\end{align}
where the second inequality is due to the traditional Fano's inequality \cite{Cover--Thomas2006}, with $U\sim \text{\rm Uniform}([M])$ and $Y\mid U\sim Q_{U}$. To proceed, we introduce a few notations: let $R_i$ be the distribution of $a(X)$ with $X\sim P_{\theta_i}$, $R$ be the distribution of $a(X)$ with $X\sim M^{-1}\sum_{i=1}^M P_{\theta_i}$, and $Q$ be the restriction of the distribution $R$ to the set $\calA_0$. Then 
\begin{align*}
I(U; Y) &\stepa{\le} \bE_U[D_{\text{KL}}(Q_U \| Q) ] \\
&\stepb{\le} \bE_U\left[\frac{1}{P_{\theta_U}(a(X)\in \calA_0)}\cdot D_{\text{KL}}(R_U \| R)  \right] \\
&\stepc{\le} \frac{1}{p_{\min}}\cdot \bE_U[D_{\text{KL}}(R_U \| R)] \\
&\stepd{=} \frac{I(U;a(X))}{p_{\min}}\\
&\stepe{\le} \frac{I(U;X)}{p_{\min}}, 
\end{align*}
where (a) is due to the variational representation of the mutual information $I(U;Y) = \min_{Q_Y} \bE_U[D_{\text{KL}}(P_{Y|U}\|Q_Y)]$, (b) follows from the data-processing property of the KL divergence $D_{\text{KL}}(P\|Q)\ge P(A)\cdot D_{\text{KL}}(P_{\cdot\mid A}\|Q_{\cdot\mid A})$, (c) is due to the assumption of Lemma \ref{lemma.fano_generalized}, (d) is the definition of the mutual information, and (e) is the data-processing property of the mutual information. Now combining the above inequality with \eqref{eq:test} completes the proof of Lemma \ref{lemma.fano_generalized}. 

\subsection{Proof of Adaptation Lower Bound}\label{subsec.detailedproof}
Recall that to formulate our adaptive property estimation problem in the general framework of \eqref{eq:minimax_adaptive}, we identify $\theta\in \Theta$ and $a\in \calA$ with the distributions $p, \widehat{p}\in \calM_k$, and $\Theta=\calA = \calM_k$. Moreover, the loss function is the absolute difference in the property value $L_F(p, \widehat{p}) = |F(p) - F(\widehat{p})|$, and the family of losses is $\calL = \{L_F: F \text{ is a 1-Lipschitz property} \}$. In this section, we apply Lemma \ref{lemma.fano_generalized} to a suitable choice of distributions $p_1,\cdots,p_M \in \calM_k$ and 1-Lipschitz properties $F_1, \cdots, F_M\in \calF_{\text{Lip}}$, and prove the target adaptation lower bound in \eqref{eq:thm1_target}. 

Without loss of generality we assume that $k=2k_0$ is an even integer. Consider the following distribution $p_0 = (p_{0,1},\cdots,p_{0,k})\in \calM_k$ serving as the ``center" of all hypotheses: 
\begin{align*}
p_0 = \left(\frac{1}{2k}, \frac{1}{2k} + \frac{1}{k(k-1)}, \frac{1}{2k}+\frac{2}{k(k-1)}, \cdots, \frac{3}{2k} \right). 
\end{align*}
Fix a parameter
\begin{align}\label{eq:delta_condition}
\delta \in \left(0, \frac{1}{4k(k-1)}\right)
\end{align}
to be chosen later, for each $u\in \calU \triangleq \{\pm 1 \}^{k_0}$ we also associate a distribution $p_u = (p_{u,1},\cdots,p_{u,k})\in \calM_k$ with
\begin{align*}
p_{u,i} = p_{0,i} + u_i\delta, \quad p_{u,k_0+i} = p_{0,k_0+i} - u_i\delta, \quad \forall i\in [k_0].
\end{align*}
Clearly each $p_u$ is a valid probability distribution, and this is known as the Paninski's construction \cite{paninski2008coincidence}. By the Gilbert--Varshamov bound, there exists $\calU_0 \subseteq \calU$ such that the minimum pairwise Hamming distance between distinct elements of $\calU_0$ is at least $k_0/5$, and $|\calU_0|\ge \exp(k_0/8)$. We will set $\{p_u\}_{u\in \calU_0}$ as the parameters $\theta_1,\cdots,\theta_M$ in Lemma \ref{lemma.fano_generalized}, with $M = |\calU_0| \ge \exp(k_0/8)$. 

For each $u\in \calU_0$, we also need to specify the associated loss, or equivalently the $1$-Lipschitz property $F_u\in \calF_{\text{Lip}}$. The detailed choice of $F_u$ is given by
\begin{align*}
F_u(p) = \sum_{i=1}^k f_u(p_i) = \sum_{i=1}^k \min_{j\in [k]} |p_i - p_{u,j}|, \quad p = (p_1,\cdots,p_k), u \in \calU_0. 
\end{align*}
As the map $x\mapsto |x-x_0|$ is $1$-Lipschitz for any $x_0\in \bR$, and the pointwise minimum of 1-Lipschitz functions is still 1-Lipschitz, each $F_u$ is a valid $1$-Lipschitz property. 

Finally, to apply Lemma \ref{lemma.fano_generalized}, it remains to specify the subset $\calA_0$. For each $i\in [k]$, let $I_i$ be the open interval $(p_{0,i} - 1/(2k(k-1)), p_{0,i} + 1/(2k(k-1)))$; clearly $I_1,\cdots,I_k$ are disjoint intervals by the definition of $p_0$. Now we define $\calA_0$ as
\begin{align*}
\calA_0 \triangleq \left\{q=(q_1,\cdots,q_k)\in \calM_k: \sum_{i=1}^k \prod_{j=1}^k \1(q_j \notin I_i) \le \frac{k}{10} \right\}. 
\end{align*} 
In other words, the subset $\calA_0$ consists of all probability vectors which intersect with at least $9/10$ of the intervals $I_1,\cdots,I_k$. 

With the above construction and definitions, we are about to use Lemma \ref{lemma.fano_generalized} for the adaptation lower bound. Specifically, we are left with three tasks: to lower bound the separation parameter $\Delta$, to lower bound the minimum probability $p_{\min}$ for all estimators $\widehat{p}\in \calP$, and to upper bound the mutual information $I(U;X^n)$. 

\vspace{1em}
\noindent\textbf{Lower bound of $\Delta$.} First, we aim to find a lower bound of $|F_u(q) - F_u(p_u)| + |F_{u'}(q) - F_{u'}(p_{u'})|$ for all $q\in\calA_0$ and $u\neq u'\in \calU_0$. By construction of $F_u$, it is clear that $F_u(p_u) = 0$ for all $u\in \calU_0$, and the above quantity can be written as
\begin{align*}
|F_u(q) - F_u(p_u)| + |F_{u'}(q) - F_{u'}(p_{u'})| = \sum_{i=1}^k \left(\min_{j\in [k]} |q_i - p_{u,j}| + \min_{j\in [k]} |q_i - p_{u',j}| \right). 
\end{align*}
One could check the following simple fact: if $q_i\in I_{j(i)}$ for some $j(i)\in [k]$, then 
\begin{align*}
\min_{j\in [k]} |q_i - p_{u,j}| + \min_{j\in [k]} |q_i - p_{u',j}| \ge |p_{u,j(i)} - p_{u', j(i)}| \in \{0, 2\delta\}. 
\end{align*}
By the definition of $q\in\calA_0$, we know that the set $\{j(i)\}_{i\in [k]}$ contains at least $9k/10$ elements of $[k]$. Moreover, by the minimum distance property of $\calU_0$, for any $u\neq u'\in \calU_0$, there are at least $k/5$ indices $j\in [k]$ such that $|p_{u,j} - p_{u', j}| = 2\delta$. By an inclusion-exclusion principle, there are at least $9k/10 + k/5 - k = k/10$ elements in the set $\{j(i)\}_{i\in [k]}$ such that $|p_{u,j(i)} - p_{u', j(i)}| = 2\delta$, and therefore
\begin{align*}
|F_u(q) - F_u(p_u)| + |F_{u'}(q) - F_{u'}(p_{u'})| \ge \frac{k}{10}\cdot 2\delta = \frac{k\delta}{5}, \quad \forall u\neq u'\in \calU_0, q\in \calA_0. 
\end{align*}
In other words, $\Delta \ge k\delta/5$ in Lemma \ref{lemma.fano_generalized}. 

\vspace{1em}
\noindent\textbf{Lower bound of $p_{\min}$.} Next, we lower bound the probability $\bP_{p_u}(\widehat{p}(X)\in \calA_0)$ for all $\widehat{p}\in \calP$ and $u\in \calU_0$. Here we need to use the definition of $\calP$ in Assumption \ref{assump}. Assume without loss of generality that $\widehat{p}_1\le \cdots \le \widehat{p}_k$, as any permutation of $\widehat{p}$ does not affect the validity of Assumption \ref{assump}. Also, by the definition of $p_u$ and the choice of $\delta$ in \eqref{eq:delta_condition}, the entries of each $p_u$ are monotonically increasing as well. Consequently, choosing $p = p_u$ in Assumption \ref{assump} gives
\begin{align*}
\bE_{p_u}\left[\sum_{i=1}^k |\widehat{p}_i - p_{u,i}| \right] \le A(n)\sqrt{\frac{k}{n}}. 
\end{align*}
On the other hand, if the event $\widehat{p}\notin \calA_0$ occurs, there are at least $k/10$ indices $i\in [k]$ such that $\widehat{p}_j \notin I_i$ for all $j\in [k]$. Consequently, for such an index $i$, one has $|\widehat{p}_i - p_{u,i}| \ge 1/(2k(k-1)) - \delta \ge 1/(4k(k-1))$ by the choice of $\delta$ in \eqref{eq:delta_condition}. Therefore, 
\begin{align*}
\sum_{i=1}^k |\widehat{p}_i - p_{u,i}| \ge \frac{k}{10}\cdot \frac{1}{4k(k-1)}\cdot \1(\widehat{p}\notin \calA_0) \ge \frac{1}{40k}\cdot \1(\widehat{p}\notin \calA_0).
\end{align*}
Combining the above two inequalities, we conclude that
\begin{align*}
\sup_{\widehat{p}\in \calP}\max_{u\in \calU_0} \bP_{p_u}(\widehat{p}(X)\notin \calA_0) \le 40A(n)\cdot \sqrt{\frac{k^3}{n}}, 
\end{align*}
which is far smaller than $1$ as $k\ll n^{1/3}$ and the assumption $A(n)\ll n^{\delta}$ for all $\delta>0$. Consequently, we may choose $p_{\min} \ge 1/2$. 

\vspace{1em}
\noindent\textbf{Upper bound of $I(U;X^n)$.} The upper bound of the mutual information could be established in a similar way as \cite{han2018local}. Specifically, the following chain of inequalities holds:
\begin{align*}
I(U;X^n) &\stepa{\le} \bE_U[D_{\text{KL}}(p_U^{\otimes n} \| p_0^{\otimes n}) ] \\
&\stepb{=} n\cdot \bE_U[D_{\text{KL}}(p_U\| p_0) ] \\
&\stepc{\le} n\cdot \bE_U\left[\sum_{i=1}^k \frac{(p_{U,i} - p_{0,i})^2}{p_{0,i}} \right] \\
&\stepd{\le} 2nk^2\delta^2, 
\end{align*} 
where (a) is due to the variational representation of the mutual information $I(U;X) = \min_{Q_X} \bE_U[D_{\text{KL}}(P_{X|U}\|Q_X)]$ and the fact that $P_{X^n|U} = p_U^{\otimes n}$, (b) follows from the chain rule of the KL divergence, (c) uses the inequality $D_{\text{KL}}(P \| Q)\le \chi^2(P\|Q)$, and (d) follows from $\min_{i\in [k]} p_{0,i}\ge 1/(2k)$ and simple algebra. Consequently, the mutual information could be upper bounded as $I(U;X^n)\le 2nk^2\delta^2$. 

\vspace{1em}
Combining the above analysis, Lemma \ref{lemma.fano_generalized} gives that
\begin{align*}
\inf_{\widehat{p}\in \calP}\sup_{F \in \calF_{\text{\rm Lip}}}\sup_{p\in \calM_k} \bE_{p} |F(\widehat{p}) - F(p)| \ge \frac{k\delta}{10}\left(\frac{1}{2} - \frac{2nk^2\delta^2+\log 2}{k/20}\right). 
\end{align*}
Consequently, choosing $\delta = c/\sqrt{nk}$ for a small enough constant $c>0$ completes the proof of the target lower bound \eqref{eq:thm1_target} (note that the condition \eqref{eq:delta_condition} on $\delta$ is also fulfilled as $k\ll n^{1/3}$). 

\section{Proof of Theorem \ref{thm:PML}}\label{sec.theorem2}
This section is devoted to the proof of Theorem \ref{thm:PML}. The proof consists of two steps: first, we show that the PML distribution belongs to the class $\calP$ in Assumption \ref{assump}, and therefore the adaptation lower bound of Theorem \ref{thm:main} holds for the PML estimator; second, we argue by contradiction that if Theorem \ref{thm:PML} is false, then the PML plug-in approach will also achieve the rate-optimal minimax rate for all $1$-Lipschitz properties for some $k\ll n^{1/3}$, a contradiction to Theorem \ref{thm:main}. 

\vspace{1em}
\noindent\textbf{Step I: show that $p^{\text{\rm PML}}\in \calP$.} First, for the empirical distribution $\widehat{p}$, \cite{han2015minimax} shows that
\begin{align*}
\sup_{p\in \calM_k} \bE_p\left[\sum_{i=1}^k |\widehat{p}_i - p_i|\right] \le \sqrt{\frac{k}{n}}. 
\end{align*}
Moreover, a single perturbation of the observations $X_1,\cdots,X_n$ only changes the quantity $\sum_{i=1}^k |\widehat{p}_i - p_i|$ by at most $2/n$. Hence, by McDiarmid's inequality, we have
\begin{align*}
\sup_{p\in \calM_k} \bP_p \left[\min_{\sigma\in \calS_k} \sum_{i=1}^k |\widehat{p}_{\sigma(i)} - p_i| \ge \sqrt{\frac{k}{n}} + \varepsilon \right] \le 2\exp\left(-\frac{n\varepsilon^2}{2}\right)
\end{align*}
for every $\varepsilon>0$. As for the PML distribution, the competitive analysis of \cite{acharya2017unified} shows that
\begin{align*}
\sup_{p\in \calM_k} \bP_p \left[\min_{\sigma\in \calS_k} \sum_{i=1}^k |p_{\sigma(i)}^{\text{PML}} - p_i| \ge 2\varepsilon \right] \le |\Phi_{n,k}|\cdot \sup_{p\in \calM_k} \bP_p \left[\min_{\sigma\in \calS_k} \sum_{i=1}^k |\widehat{p}_{\sigma(i)} - p_i| \ge \varepsilon \right], 
\end{align*}
where $|\Phi_{n,k}|$ is the cardinality of all possible profiles with length $n$ and support size $k$. Note that trivially $|\Phi_{n,k}|\le (n+1)^k$ holds, the above two inequalities lead to
\begin{align}\label{eq:PML_sorted_dist}
\sup_{p\in \calM_k} \bP_p \left[\min_{\sigma\in \calS_k} \sum_{i=1}^k |p_{\sigma(i)}^{\text{PML}} - p_i| \ge 2\varepsilon \right]\le \min\left\{1,2\exp\left(k\log(n+1) - \frac{n}{2}\left(\varepsilon - \sqrt{\frac{k}{n}}\right)_+^2 \right)\right\}. 
\end{align}
Now integrating the RHS of \eqref{eq:PML_sorted_dist} over $\varepsilon\in (0,\infty)$ gives that $p^{\text{PML}}\in \calP$ with $A(n) = O(\sqrt{\log n})$. 

\vspace{1em}
\noindent\textbf{Step II: proof by contradiction.} Assume by contradiction that Theorem \ref{thm:PML} is false, i.e. there exists an absolute constant $c_0$ such that for some large enough $n$, it holds that
\begin{align}\label{eq:contradiction}
\sup_{p\in \calM_k} \bP_p( |F(p^{\text{\rm PML}})- F(p)| \ge C\varepsilon ) \le \exp(c_0n^{1/3-c_1})\cdot \left(\inf_{\widehat{F}}\sup_{p\in \calM_k} \bP_p( |\widehat{F} - F(p)| \ge \varepsilon )\right)^{1-c_2}
\end{align}
for all $k\in \NN, \varepsilon>0$ and $F\in \calF_{\text{Lip}}$. For any $\varepsilon\gg n^{-1/2}$ and $k\gg 1$, it was shown in \cite{hao2019unified} that the minimax error probability for any $1$-Lipschitz property estimation is at most
\begin{align*}
\inf_{\widehat{F}}\sup_{p\in \calM_k} \bP_p( |\widehat{F} - F(p)| \ge \varepsilon ) \le 2\exp\left(-c_\delta n^{1-\delta}\left(\varepsilon - d_\delta\sqrt{\frac{k}{n\log n}} \right)_+^2 \right), 
\end{align*}
for an arbitrary constant $\delta>0$ and constants $c_\delta, d_\delta>0$ depending only on $\delta$. Consequently, \eqref{eq:contradiction} implies that 
\begin{align*}
\sup_{F\in\calF_{\text{Lip}}}\sup_{p\in \calM_k} \bP_p( |F(p^{\text{\rm PML}})- F(p)| \ge C\varepsilon ) \le 2\exp\left(c_0n^{1/3-c_1} - (1-c_2)c_\delta n^{1-\delta}\left(\varepsilon - d_\delta\sqrt{\frac{k}{n\log n}} \right)_+^2 \right). 
\end{align*} 
Choosing $\delta<c_1/4$, $\varepsilon = 2d_\delta\sqrt{k/(n\log n)}$ and $k \asymp n^{1/3-c_1/2}$, the above inequality shows that there exists an absolute constant $c_0'>0$ depending only on $(c_0,c_1,c_2,C)$ such that 
\begin{align*}
\sup_{F\in\calF_{\text{Lip}}}\sup_{p\in \calM_k} \bP_p\left( |F(p^{\text{\rm PML}})- F(p)| \ge \frac{1}{c_0'}\sqrt{\frac{k}{n\log n}} \right) \le 2\exp\left(-c_0'n^{1/3-c_1}\right). 
\end{align*}
Hence, using that $\bE|X| \le t + \|X\|_\infty\cdot \bP(|X|\ge t)$ for any $t>0$ implies that for $k\asymp n^{1/3 - c_1/2}$ and $n$ tending to infinity (possibly along some subsequence), we arrive at
\begin{align*}
\sup_{F\in\calF_{\text{Lip}}}\sup_{p\in \calM_k} \bE_p|F(p^{\text{\rm PML}})- F(p)| \lesssim \sqrt{\frac{k}{n\log n}}, 
\end{align*} 
a contradiction to Theorem \ref{thm:main} as $p^{\text{PML}}\in \calP$. Therefore, the inequality \eqref{eq:contradiction} does not hold, and the proof of Theorem \ref{thm:PML} is completed. 

\section{Conclusion and Open Problems}\label{sec.conclusion}
In this paper we showed that there is a high-accuracy limitation for general adaptive approaches of property estimation, which in turn implied tight lower bounds for the known adaptive approaches such as the PML and LMM. A number of directions could be of interest. First, we believe that Assumption \ref{assump} is an artifact of our proof and unnecessary for Theorem \ref{thm:main} to hold, and a better choice of the loss functions in Lemma \ref{lemma.fano_generalized} could remove this assumption. Second, the adaptation lower bound for PML does not rule out the possibility that PML could be fully optimal for \emph{certain properties}. However, to show this, one need to go beyond the competitive analysis of the PML and seek for additional properties. Third, our current lower bound for PML only shows the existence of a property requiring $\varepsilon\gg n^{-1/3}$ for the PML to be optimal, and it is interesting to construct such a property explicitly. 

\section*{Acknowledgement} 
Yanjun Han is grateful to Kirankumar Shiragur for helpful discussions, and anonymous reviewers for their valuable comments to improve the presentation of this paper.

\bibliographystyle{alpha}
\bibliography{di}

\end{document}